\documentclass[12pt,a4paper]{article}
\usepackage{amsmath,amsthm,amssymb,latexsym,epic}
\newcommand{\IS}{{\mathcal{I}}}
\renewcommand{\S}{{\mathcal{S}}}
\newcommand{\dom}{{\mathrm{dom}}}
\newcommand{\rank}{{\mathrm{rank}}}
\newcommand{\im}{{\mathrm{im}}}
\newcommand{\SIS}{\mathcal{I}_n\setminus \mathcal{S}}
\newtheorem{theorem}{Theorem}
\newtheorem{lemma}{Lemma}

\author{Karin Cvetko-Vah, Damjana Kokol Bukov\v{s}ek, Toma\v{z} Ko\v{s}ir\\ and Ganna Kudryavtseva}
\title{Semitransitive subsemigroups of the singular part of the finite symmetric inverse semigroup}
\begin{document}
\maketitle
\begin{abstract}
We prove that the minimal cardinality of the semitransitive subsemigroup in the singular part $\IS_n\setminus \S_n$ of the symmetric inverse semigroup $\IS_n$ is $2n-p+1$, where $p$ is the greatest proper divisor of $n$, and classify all semitransitive subsemigroups of this minimal cardinality.
\end{abstract}

\section{Introduction}

A semigroup $S$ of transformations of the set $X$ is called {\em semitransitive} \cite{RT} if for every $x,y\in X$ there is $\varphi\in S$ such that either $x\varphi=y$ or $y\varphi=x$.
A number of papers has been published which are devoted to the study of semitransitive actions of spaces of linear operators \cite{BDKKO, BDKKOR, BGMRT, RaT, RT, Bled}. In
\cite{C-O} we initiated the study of semitransitive semigroups of transformations.

In this paper we continue the research started in \cite{C-O}. There we described the semitransitive subsemigroups of $\IS_n$ of the minimal cardinality. We proved that every semitransitive subsemigroup of $\IS_n$ of the minimal cardinality contains the identity (which is the only idempotent), and so is not contained in the singular part of $\IS_n$. By the {\em singular part} of the finite symmetric  inverse semigroup $\IS_n$ we mean the semigroup $\IS_n\setminus \S_n$, that is the maximal ideal of $\IS_n$ consisting of all non-invertible elements of $\IS_n$. This observation shows that a straightforward reduction of our problem to the corresponding result for $\IS_n$ is not possible. This is not very  surprising because  singular parts of transformation semigroups have richer combinatorics, see for example, \cite{East, MM}.

The structure of the paper is as follows. In Section \ref{s2} we study semitransitive subsemigroup of $\IS_n\setminus \S_n$ of small cardinality ($\leq 2n$) and show that any such a semigroup must contain exactly two non-zero idempotents, and each transitivity block is a subset of the domain of one of these idempotents. In Section \ref{s3} we give a lower bound for the cardinality of the set of nilpotent elements, and find some restrictions to the cardinalities of the transitivity blocks. This allows us to prove that the cardinality of a semitransitive subsemigroup in $\IS_n\setminus \S_n$ can not be less than $2n-p+1$, where $p$ is the greatest proper divisor of $n$. Finally, in Section \ref{s4}, we construct five types of semitransitive subsemigroups in $\IS_n\setminus \S_n$ of the cardinality $2n-p+1$, showing that the lower bound above is sharp. We prove that each semitransitive subsemigroups in $\IS_n\setminus \S_n$ of the minimal cardinality is of one of the five types constructed.

\section{Idempotents and their position with respect to transitivity blocks}\label{s2}

Fix $S$ to be a semitransitive semigroup contained in $\SIS_n$ and $|S|\leq 2n$. From the description of transitive subsemigroups of $\IS_n$ (\cite[Theorem 2.5, Remark 2.7]{C-O}, this result is recalled in detail in Section \ref{s4} below) it follows that $S$ cannot be transitive, and therefore by \cite[Theorem 3.1]{C-O} $S$ contains the zero, which we denote by $0$.

Let $X$ denote the underlying set on which $\IS_n$ acts.

For $x,y\in X$ set $x\geq y$ provided that there is $\varphi\in S$ such that $x\varphi=y$. The relation $\geq$ is a linear preorder which we call the {\em semitransitivity preorder} induced by $S$. For each possible $s$ define $X_s$ to be the set of all $x\in X$ such that for $y\in X$ there is $\varphi\in S$ with $y\varphi=x$ if and only if $y\in X_1\cup\dots\cup X_s$, and $\varphi\in S$ with $x\varphi=y$ if and only if  $y\in X_s\cup\dots\cup X_m$.  We call the blocks $X_i$ {\em the transitivity blocks} (or just blocks).
 Let $X_1\geq X_2\geq \cdots\geq X_m$ of $X$ be the ordering of blocks induced by semitransitivity preorder on $X$.

\begin{lemma} \label{lem:two_idemp}$S$ contains exactly two non-zero idempotents. Their domains are disjoint, and the union of their domains is $X$.
\end{lemma}

\begin{proof}
Since for each $x\in X$ there is an element in $S$ which maps $x$ to itself, it follows that there is an idempotent $e_x\in S$ such that $x\in \dom(e_x)$. So if $S$ had only one idempotent it would be the identity map on $X$, which is not contained in the $\SIS_n$. So $S$ contains at least two idempotents, and the union of all domains of idempotents in $S$ is $X$. Now it is enough to prove that $S$ has exactly two non-zero idempotents (then their domains are automatically disjoint as otherwise their product would be a third non-zero idempotent).

Let  $E=E(S)$ be the set of idempotents of $S$  ($E$ is naturally partially ordered with the relation $e\leq f$ if and only if $ef=fe=e$) and $e_1,\dots, e_t$ the non-zero elements in $E$.   Assume that $t\geq 3$. For $e\in E$ set $M_e=\dom(e)$.  For $e\in E\setminus\{0\}$ set $$M_e'=M_e\setminus \cup_{f<e}M_f, \,\, m_e=|M_e'|.$$

For $e,f\in E\setminus\{0\}$ let
$$
S_{e,f}=eSf\setminus \bigl((\cup_{a<e}aSf)\cup (\cup_{b<f}eSb)\bigr).
$$

By the definition, $S_{e,f}$ is either empty or it is the set of all elements $\varphi$ in $S$ such that $\dom(\varphi)\subseteq M_e$, $\im(\varphi)\subseteq M_f$ and $\dom(\varphi)\cap M_e'\neq \varnothing$ and $\im(\varphi)\cap M_f'\neq \varnothing$. Observe that $S_{e,f}\cup S_{f,e}\neq \varnothing$ and
the sets $S_{e,f}\cup S_{f,e}$ and $S_{g,h}\cup S_{h,g}$ are disjoint unless $\{e,f\}=\{g,h\}$.

Let $e\neq f$ be two nonzero idempotents. Let us estimate the cardinality of $S_{e,f}\cup S_{f,e}$. Assume that $m_e\geq m_f$. Fix $x\in M_f'$ and for each $y\in M_e'$ let $\varphi_{x,y}\in S_{e,f}\cup S_{f,e}$ be such that $x\varphi_{x,y}=y$ or $y\varphi_{x,y}=x$. Suppose $\varphi_{x,y'}=\varphi_{x,y''}$ for $y',y''\in M_e'$. Then $\dom(\varphi_{x,y'})$ contains $x$ and at least one of $y',y''$. Assume that $\dom(\varphi_{x,y'})\supseteq \{x,y'\}$. Then $\im(\varphi_{x,y'})\supseteq \{x,y''\}$. But this contradicts the fact that we have either $\dom(\varphi_{x,y'})\subseteq M_e$ and $\im(\varphi_{x,y'})\subseteq M_f$ or $\dom(\varphi_{x,y'})\subseteq M_f$ and $\im(\varphi_{x,y'})\subseteq M_e$. It follows that all the elements $\varphi_{x,y}$, $y\in M_e'$,  are pairwise different and thus $|S_{e,f}\cup S_{f,e}|\geq m_e$.

Denote the cardinalities of the sets $M'_e$, $e\in E\setminus\{0\}$, by $m_1,\dots, m_t$. Suppose that $m_t\geq \cdots \geq m_1$. We have that  $m_1+\cdots +m_t=n$.  Bearing in mind that $t\geq 3$ we calculate that  the union of all the sets of the form $S_{e,f}\cup S_{f,e}$, $e\neq f$, contains at least
$$
(t-1)m_t+(t-2)m_{t-1}+\cdots +m_2\geq m_t+\cdots +m_1=n
$$
elements.

Finally, let us estimate the cardinalities of the sets $S_{e,e}$. Let $M'_e=\{x_1, \dots, x_{m_e}\}$. Let $x\in M'_e$ be such that $x\geq x'$ for any $x'\in M'_e$. It follows that for any $x'\in M'_e$ there is $\varphi_{x,x'} \in S_{e,e}$ such that $x\varphi=x'$. Since all the elements $\varphi_{x,x'}$ are pairwise different, we have that $|S_{e,e}|\geq m_e$ and the sum of cardinalities of all  $S_{e,e}$'s is at least $n$.

It follows, that $S$ must contain at least $n+n=2n$ non-zero elements. This completes the proof.
\end{proof}

Denote by $g$ and $h$ the two non-zero idempotents of $S$.

Let $\alpha\in S$. We will say that $\alpha$ {\em has an arrow} $x\to y$ provided that $x\in\dom(\alpha)$ and $x\alpha=y$. The notation stems from considering the graph of action of $\alpha$: $\alpha$ has the arrow $x\to y$ whenever the graph of action of $\alpha$ has a directed edge $(x,y)$. In the case when $x\in X_p$ and $y\in X_q$ we will say that the arrow $x\to y$ of $\alpha$ is from the block $X_p$ to the block $X_q$.

\begin{lemma}\label{lem:aux} For any $x,y$ belonging to the same block there is $\alpha\in S$ which has an arrow $x\to y$ and does not have arrows from $X_p$ to $X_q$ with $p\neq q$.
\end{lemma}

\begin{proof} Suppose $\alpha$ has an arrow $x\to y$ where $x,y\in X_r$ and an arrow $u\to v$ with $u\in X_p$, $v\in X_q$, $p>q$. There is $\beta$  such that it has an arrow $y\to x$. Let $t$ be such that the element $e=(\alpha\beta)^t$ is an idempotent. Then $e$ does not contain $u$ in its domain. Therefore, $e\alpha$, while containing an arrow $x\to y$, does not contain an arrow $u\to v$. \end{proof}

In the following lemma we establish the connection between the  transitivity blocks of $S$ and the domains of the two non-zero idempotents $g$ and $h$.

\begin{lemma}\label{lemma:lem2} Let $e\in S$ be an idempotent. Then for any block $X_s$ we have that either $\dom(e)\cap X_s=\varnothing$ or $\dom(e)\cap X_s=X_s$.
\end{lemma}

\begin{proof}
Apply induction on $s$. Suppose $s=1$.  If $|X_1|=1$, the statement holds by a trivial argument. Let $|X_1|\geq2$. Assume the converse to the claim for $X_1$. Consider $x,y\in X_1$, $x\neq y$, such that  $x\in \dom(g)$ and $y\in\dom(h)$. For each $t\in X$ there are $\gamma_{x,t}$, $\delta_{y,t}$ with $x\gamma_{x,t}=t$ and $y\delta_{y,t}=t$. The elements $g\gamma_{x,t}$, $h\delta_{y,t}$, $t\in X$, are pairwise distinct. It follows that $S$ has at least $2n$ non-zero elements. A contradiction. It follows that either $\dom(g)\supseteq X_1$ or $\dom(h)\supseteq X_1$.

Let $s\geq 2$. Again, if  $|X_s|=1$, there is nothing to prove. Suppose $|X_s|\geq 2$ and assume the converse to the claim for $X_s$. Let $x,y\in X_s$, $x\neq y$, be such that  $x\in \dom(g)$ and $y\in\dom(h)$.  For $t\in X_1\cup\dots\cup X_{s-1}$ let $\gamma_{t,x}$, $\gamma_{t,y}$ be such that $t\gamma_{t,x}=x$, $t\gamma_{t,y}=y$ and for $z\in X_s\cup\dots\cup X_m$ let $\delta_{x,z}$, $\delta_{y,z}$ be such that $x\delta_{x,z}=z$ and $y\delta_{y,z}=z$. Choose $\alpha$ and $\beta$ so that
$y\alpha=x$, $x\beta=y$ and $\alpha$ and $\beta$ do not have arrows between different blocks (this is possible to do by Lemma~\ref{lem:aux}). Then the element
$\delta'_{y,z}=h\alpha g\beta \delta_{y,z}$
has an arrow $y\to z$ and, in view of the inductive hypothesis, is such that $\dom(\delta'_{y,z})\cap (X_1\cup\dots\cup X_{s-1})=\varnothing$. Similarly we construct the elements $\delta'_{x,z}$ such that $\delta'_{x,z}$ has an arrow $x\to z$ and $\dom(\delta'_{x,z})\cap  (X_1\cup\dots\cup X_{s-1})=\varnothing$.
Then the $2n$ elements $\delta'_{x,z}$, $\delta'_{y,z}$, $z\in X_s\cup\dots\cup X_m$, and $\gamma_{t,x}g$, $\gamma_{t,y}h$, $t\in X_1\cup\dots\cup X_{s-1}$, are pairwise distinct. We again obtain that $|S|$ has at least $2n$ non-zero elements. This contradiction completes the proof.
\end{proof}

\section{Nilpotent elements and the lower bound for the cardinality}\label{s3}

\begin{lemma} \label{lemma:lem3} For any non-zero idempotent $e\in S$ and $\varphi\in eSe$ we have that $\varphi$ either belongs to the group of units of $eSe$ or is nilpotent.
\end{lemma}

\begin{proof} Let $\varphi\in eSe$ and $e_{\varphi}$ be its idempotent power. Since $\dom(e_{\varphi})\subseteq \dom(e)$ then by Lemma \ref{lem:two_idemp} we have either $e_{\varphi}=0$ or $e_{\varphi}=e$. The statement follows.
\end{proof}

$S$ has a subsemigroup $S'=gSg\cup gSh\cup hSg \cup hSh$ which is also semitransitive and has the same semitransitivity preorder and transitivity blocks as $S$. Because the cardinality of $S'$ is not greater than that of $S$, and the limitation to the cardinality was the only one imposed on $S$, all the statements we proved above for $S$, hold also for $S'$.

Denote by $N_S$  the set of nilpotent elements of $S'$.

\begin{lemma}\label{lem:new}
Let $\varphi\in N_S$. Then for any $x\in\dom\varphi$ we have that if $x\in X_i$ and $x\varphi\in X_j$ then $j>i$. In other words, a nilpotent element of $S'$ does not have arrows from a block to itself.
\end{lemma}

\begin{proof} Suppose that there is $z\in\dom(\varphi)$ such that $z\in\dom(\varphi)\cap X_i$ and $s=z\varphi\in X_i$.  Consider $\psi$ such that $s\psi=z$. Then $z$ belongs to the domain of the idempotent power of $\varphi\psi$. It follows that the idempotent power of $\varphi\psi$ is a non-zero idempotent, say, $e$. It follows that $\dom(e)=\dom(\varphi)$. By a similar argument, the idempotent power of $\psi\varphi$ also equals $e$, and then $\im(e)=\im(\varphi)$. It follows that $\varphi$ acts bijectively on the set $\dom(e)$, which implies that some its power is $e$. A contradiction.
\end{proof}

Let  $N_{1,2}=gSh$, $N_{2,1}=hSg$
and $N= \left(N_{1,2}\cup N_{2,1}\right)\setminus \{0\}$. Observe that all elements of $N$ are nilpotent and  if $\varphi\in N$ and $x\in \dom(\varphi)\cap X_i$ then $x\varphi\in X_j$ with $j>i$. Let $t_i=|X_i|, 1\leq i\leq m$. Set $t=\min_{1\leq i\leq m}t_i$.

\begin{lemma}\label{lemma:lem4}
$|N|\geq n-t$.
\end{lemma}
\begin{proof}
Let $A=\{i:\, X_i\subseteq \dom(g)\}$, $B=A^c=\{i:\, X_i\subseteq \dom(h)\}$. For each $r$,
$1\leq r\leq m-1$, we consider the sets
$$
M^r_{1,2}=N_{1,2}\cap(S_{\bf n}^r\setminus S_{\bf n}^{r+1}),
M^r_{2,1}=N_{2,1}\cap(S_{\bf n}^r\setminus S_{\bf n}^{r+1}),
M^r=M^r_{1,2}\cup M^r_{2,1}.
$$
 Set also $M^m=\varnothing$. Observe that $M^r_{1,2}\cap M^r_{2,1}=\varnothing$ for each $r$, $M^i\cap M^j=\varnothing$ for $i\neq j$ and
$N= M^1\cup\dots\cup M^m$.

The construction implies that for $x\in X_i$, $y\in X_j$, $j-i=l\geq  1$, where $i\in A$, $j\in B$ or  $i\in B$, $j\in A$, there is an element  $\alpha \in M^l$ such that
$x\alpha=y$.

For every $i\in A$ define the set
$B_i=\{j-i: \, j\in B\}$. Observe that for $j\in B_i$ we have the estimates
\begin{gather}
|M^{j}_{1,2}|\geq\max\{t_i,t_{i+j}\}, \text{ if }
j>0,\label{g1}\\
|M^{-j}_{2,1}|\geq\max\{t_i,t_{i+j}\}, \text{ if } j<0.\label{g2}
\end{gather}

We claim that the union of all sets $B_i$, $i\in A$, contains at least $|B|+|A|-1=m-1$ elements. Indeed, if $i_1< i_2< \dots<  i_{|A|}$ are all the elements of $A$, then $B_{i_1}$ contributes to the union $|B|$ elements. Then each of $B_{i_k}$, $k\geq 2$, contributes at least one new element as the minimum element of $B_{i_k}$, $k\geq 2$, is smaller that any element from the union of $B_i$'s with $i<i_k$.

Since $N$ is a disjoint union of at least $m-1$ different non-empty sets of the form $M_{2,1}^{j}$ and $M_{1,2}^{j}$ it follows from~\eqref{g1},~\eqref{g2} and $t_1+\cdots +t_m=n$ that $|N|\geq n-t$.
\end{proof}

\begin{lemma} \label{lemma:lem5} For each $i$, $1\leq i\leq m$,  $t_i$ is divisible by $t$.
\end{lemma}
\begin{proof}
Let $s$ be such that $t=t_s$. We can suppose that $s\in A$. Assume that $t_i$ is not divisible by
$t$.

{\it Case 1.} Suppose that $i\in A$. Let $n_1=\rank(g)$ and
$n_2=\rank(h)$. From~\cite[Proposiion 3.1, Theorem 3.3]{C-O} we know that
$|gSg\setminus\{0\}|\geq n_1$, $|hSh\setminus\{0\}|\geq n_2$. It is enough to show that in this case the cardinality of
$gSg$ is bigger than $n_1+t$, which, in view of
Lemma~\ref{lemma:lem4}, would lead to $|S\setminus\{0\}|\geq 2n$, which is not possible.

Denote $T=gSg$. By Lemma~\ref{lemma:lem3} we have that $T=T_{\bf g}\cup T_{\bf n}$, where
$T_{\bf g}$ is the set of all group elements, and $T_{\bf n}$ the
set of all nilpotent elements of $T$. In turn, $T_{\bf n}$ is the
disjoint union of the sets $T_j=T_{\bf n}^j\setminus
T_{\bf n}^{j+1}$, $1\leq j\leq |A|-1$. Since $|T_{\bf g}|\geq t$ it is enough to show that  $|T_{\bf n}|\geq n_1$.

Observe that for any $x\in X_{i_k}$, $y\in X_{i_{k+j}}$ there is $\varphi\in T_j$ such that $x\varphi=y$. It follows that
$$
|T_j|\geq \left\lbrace \begin{array}{ll} max\{t_{i_1}, \dots, t_{i_{|A|}}\},& j\leq\frac{|A|-1}{2};\\
max\{t_{i_1},\dots, t_{i_{|A|-j}}, t_{i_{j+1}}, \dots, t_{i_{|A|}}\},&
j>\frac{|A|-1}{2}.\end{array}\right.
$$

By adding up the cardinalities, we obtain that
\begin{equation}\label{eq:card1}
|T_{\bf n}|\geq n_1-t.
\end{equation}

Let $\varphi\in T_{\bf n}$ be such that it has an arrow from $X_i$
to $X_s$ (we assume $i<s$, the other case is treated similarly). Let $j$ be maximal such that $\varphi$ can be chosen in $T_j$ ($j=b-a$ where $X_i=X_{i_a}$, $X_s=X_{i_b}$ by the definition of $T_j$).

Let $Z=\{y\in \im(\varphi)\cap X_s: y\varphi^{-1} \in X_i\}$, $a=|Z|$. Consider first the case when $a=t$.
Let $l\geq 1$ be such that $lt<t_i<(l+1)t$. Let
$A_1=Z\varphi^{-1}$. As the group  $T_{\bf g}$ acts
transitively on $X_i$, there are $A_1,\dots, A_r\subset
X_i$ such that $r\geq l+1$, $|A_k|=a$ for all $k$'s,
$X_i=A_1\cup\dots\cup A_r$ and, in addition, for every $\psi\in T_{\bf g}$ and every possible $k$ we have that
$\psi(A_k)=A_h$ for some $h$. There are some $A_u$ and $A_v$
whose intersection is not empty. Suppose $x\in A_u\cap A_v$. Then
$x\psi\in\psi(A_u)\cap\psi(A_v)$ for every $\psi\in T_{\bf g}$. Since the group  $T_{\bf g}$ acts transitively on
blocks, it follows that every element in $X_i$ belongs to some
$A_u\cap A_v$. Similarly we show that there is $q\geq 2$ such that
every element of $X_i$ belongs to intersection of exactly $q$ of
$A_k$'s while intersection of any $(q+1)$ of $A_k$'s is empty. Then
we have
$$
ar=qt_i>qla,
$$
which implies  $r>ql\geq 2l$  and $r\geq 2l+1\geq l+2$.

Fix some $x\in X_s$ and let $g_1$, $\dots, g_s$ be the elements of $T_{\bf g}$ which map $x$ to different elements of $X_s$. Then all the elements $\varphi g_i$ are pairwise distinct, belong to $T_j$ and have the property that $X_s(\varphi g)^{-1}=A_1$. Multiplying the elements $\varphi g_i$ with different elements in $T_{\bf g}$ we obtain elements $\psi$ such that $X_s\psi^{-1}=A_k$ for all $k$. As we have at least $l+2$ different $A_k$'s it follows that  $|T_{j}|\geq (l+2)a\geq t+t_i$ which in view of~\eqref{eq:card1} implies that
$|T_{\bf n}|\geq n_1$. A contradiction. Therefore, this case is
impossible.

Consider the case $t>a$.  Then we can find $f\in T_{\bf g}$ such that $\im(\varphi f)\cap X_s\neq Z$. Multiplying the elements $\varphi$ and $\varphi f$ with different $g\in T_{\bf g}$ from the left we obtain at least $t_i$ pairwise different elements $g\varphi$ and $t_i$ pairwise different elements $g\varphi f$. It follows that $|T_{j}|\geq 2t_i\geq t+t_i$. A contradiction. Hence this case is also impossible.

{\it Case 2.} Suppose that $i\in B$. In this case we apply similar
arguments as in the previous case. In view of Lemma~\ref{lemma:lem4}
we make a conclusion that $|N|\geq n$. Then we have that $S$ has  at
least $n+n_1+n_2=2n$ non-zero elements. A contradiction.
\end{proof}

\begin{theorem}\label{th}
Let $S$ be a semitransitive semigroup contained in $\IS_n\setminus \S_n$ and $p$ be the greatest proper divisor of $n$. Then $|S|\geq 2n-p+1$.
\end{theorem}

\begin{proof}
The statement follows from Lemmas~\ref{lemma:lem4} and~\ref{lemma:lem5} and the fact that $gSg\cup hSh$ contains at least $n$ non-zero elements.
\end{proof}



\section{Classification of semitransitive subsemigroups of $\SIS_n$ of the minimal cardinality}\label{s4}

{\bf Construction.}
Let $Z$ be a finite set, $G$ a transitive permutation group of $Z$ of cardinality $|Z|$, $l\geq 2$ and $T$ a subsemigroup of $\IS_l$.  The semigroup $G\times T$ acts on the set $Z\times \{1,2\dots, l\}$ by partial permutations as follows:
$$
(\alpha,\beta)\cdot(z,i)=\left\lbrace\begin{array}{ll}(\alpha z,\beta i), & {\text{ if }}   i\in dom(\beta)\\
\text{ not defined, }& {\text { otherwise.}}\end{array}\right.
$$

This action is non faithful in general as in case when $T$ has the zero element $0$ all the elements $(\alpha, 0)$, $\alpha\in G$, act the same way: they are nowhere defined. However, the induced action of the Rees factor semigroup $(G\times T)/I$, where the ideal $I$ consists of all the elements $(\alpha, 0)$, $\alpha\in G$, is faithful and we can consider $(G\times T)/I$ as a semigroup of partial permutation of $Z\times \{1,2\dots, l\}$. We identify $Z$ with $Z_1$ and $Z\times \{1,2\dots, l\}$ with the union $Z_1\cup Z_2\cup \dots \cup Z_l$ of pairwise disjoint sets $Z_i$, $1\leq i\leq l$, each of which has cardinality $|Z|$ by assigning to $(z,i)$ the image element $z$ in the block $Z_i$ under some fixed bijection $Z\to Z_i$. Under this identification we obtain the faithful action of $(G\times T)/I$ (or of $G\times T$ if $T$ does not have the zero) on the set $Z_1\cup Z_2\cup \dots \cup Z_l$. Let $n=|Z|\cdot l$. In the sequel we will always consider the case when the semigroup $T$ does have the zero. In this case the cardinality of $(G\times T)/I$ is $|G|\cdot (|T|-1)+1$.

For a semigroup $S$ we denote by $S^1$ the semigroup $S$ with the adjoint identity element, if $S$ does not have the identity, and we write $S^1=S$, if $S$ has the identity. In \cite{C-O} the following result was proved (for the chain-cycle notation for the elements of $\IS_n$ we refer the reader to \cite{GM}):

\begin{theorem}\label{th_from_c_o}
The minimal cardinality of a semitransitive, but not transitive, subsemigroup of $\IS_n$ is $n+1$. Any such a semigroup is similar to a semigroup $(G\times T^1)/I$ with the action described above for some decomposition $Z_1\cup Z_2\cup \dots \cup Z_l$ of $\{1,2\dots, n\}$, some transitive permutation group on $Z_1$ of cardinality $|Z_1|$ and the semigroup $T$ generated by the chain $(1,2,\dots, l]$.
\end{theorem}

Let $p$ be a proper divisor of $n$. Set $m=\frac{n}{p}$. We apply the  construction above to some specific semigroups $T$ to obtain subsemigroups of $\SIS_n$ of cardinality $2n-p+1$.

{\it \bf Type 1.}  Let $X_1\subseteq X$, $|X_1|=p$, $G$ be a transitive permutation group on $X_1$ of cardinality $p$ and $T$ the subsemigroup of $\IS_m$ generated by $\varphi=(1,2](3]\dots (m]$, $\psi=(1](2,3,\dots, m]$, $g=(1)(2](3]\dots (m]$,  and $h=(1](2)\dots (m)$.  Consider some decomposition $X=X_1\cup\dots\cup X_m$.  We have $|T|=2m$ and therefore $|(G\times T)/I|=2n-p+1$.

{\it Example 1.} Let $n=8$, $p=2$, $X_1=\{1,2\}$, $X_2=\{3,4\}$,
$X_3=\{5,6\}$, $X_4=\{7,8\}$. $G=\{e, (1,2)\}$. $(G\times T)/I$ has the following $15$ elements:

\begin{equation*}
\begin{array}{ll}
\text{\bf the elements of } gSg\setminus\{0\}:&\\
(1)(2)(3](4](5](6](7](8], &  (1,2)(3](4](5](6](7](8],  \\
\text{\bf the elements of } hSh\setminus\{0\}:&\\
(1](2](3)(4)(5)(6)(7)(8), &  (1](2](3,5,7](4,6,8],  \\
(1](2](3,7](4,8](5](6], &  (1](2](3,4)(5,6)(7,8),  \\
(1](2](3,6,7](4,5,8], & (1](2](3,8](4,7](5](6],  \\
\text{\bf the elements of } gSh\setminus\{0\}:&\\
(1,3](2,4](5](6](7](8],  & (1,5](2,6](3](4](7](8], \\
(1,7](2,8](3](4](5](6],  & (1,4](2,3](5](6](7](8], \\
(1,6](2,7](3](4](7](8],  & (1,8](2,7](3](4](5](6] \\
\end{array}
\end{equation*}
and the zero.

{\it \bf Type 2.}
 Let $X_1\subseteq X$, $|X_1|=p$ and $G$ be a transitive permutation group on $X_1$ of cardinality $p$ and $T$ the subsemigroup of $\IS_m$ generated by $(1,2,\dots, m-1] (m]$, $(1](2]\dots (m-1](m)$, $(1]\dots (m-2](m-1,m]$ and $(1)(2)\dots (m-1) (m]$.  We again have $|T|=2m$ and therefore $|(G\times T)/I|=2n-p+1$. Actually, such a semigroup can be obtained from one constructed in Type 1 applying the map $\alpha\mapsto \alpha^{-1}$.

{\it \bf Type 3.} Let $X_1\subseteq X$, $|X_1|=p$ and $G$ be a transitive permutation group on $X_1$ of cardinality $p$ and $T$ the subsemigroup of $\IS_m$ generated by
$$\varphi=\left\lbrace\begin{array}{ll}(1,2](3,4]\dots (m-1,m], & \text{ if } $m$ \text{ is even},\\
(1,2](3,4]\dots (m-2,m-1](m], & \text{ if } $m$ \text{ is odd},
\end{array}\right.$$

$$\psi=\left\lbrace\begin{array}{ll}(1](2,3](4,5]\dots (m-2,m-1](m], & \text{ if } $m$ \text{ is even},\\
(1](2,3](4,5]\dots (m-1,m], & \text{ if } $m$ \text{ is odd},
\end{array}\right.$$

$$g=\left\lbrace\begin{array}{ll}(1)(2](3)(4]\dots (m-1)(m], & \text{ if } $m$ \text{ is even},\\
(1)(2](3](4]\dots (m-1](m), & \text{ if } $m$ \text{ is odd},
\end{array}\right.$$

$$h=\left\lbrace\begin{array}{ll}(1](2)(3](4)\dots (m-1](m), & \text{ if } $m$ \text{ is even},\\
(1](2)(3](4)\dots (m-1)(m], & \text{ if } $m$ \text{ is odd}.
\end{array}\right.$$

The elements of $T$ are the zero, the elements $g$, $h$ and the elements $\varphi$, $\psi$, $\varphi\psi$, $\psi\varphi$,  $\varphi\psi\varphi$, $\psi\varphi\psi$, and so on. The longest non-zero product starting with $\varphi$ has $m-1$ factors, and the longest non-zero product starting with $\psi$ has $m-2$ factors. It follows that $|T|=2m$ and hence $|(G\times T)/I|=2n-p+1$.

{\it Example 3.} Let $n=8$, $p=2$, $X_1=\{1,2\}$, $X_2=\{3,4\}$,
$X_3=\{5,6\}$, $X_4=\{7,8\}$, $G=\{e, (1,2)\}$.

\begin{equation*}
\begin{array}{ll}
\text{\bf The elements of } gSg\setminus\{0\}:&\\
(1)(2)(5)(6)(3](4](7](8], & (1,2)(5,6)(3](4](7](8], \\
(1,5](2,6](3](4](7](8], & (1,6](2,5](3](4](7](8], \\
\text{\bf the elements of } hSh\setminus\{0\}:&\\
(3)(4)(7)(8)(1](2](5](6], & (3,4)(7,8)(1](2](5](6],\\
(3,7](4,8](1](2](5](6],& (3,8](4,7](1](2](5](6], \\
\text{\bf the elements of } gSh\setminus\{0\}:&\\
(1,3](2,4](5,7](6,8], & (1,4](2,3](5,8](6,7],\\
(1,7](2,8](2](3](4](5], & (1,8](2,7](2](3](4](5], \\
\text{\bf the elements of } hSg\setminus\{0\}:&\\
(3,5](4,6](1](2](7](8], & (3,6](4,5](1](2](7](8]. \\
\end{array}
\end{equation*}

{\it \bf Type 4.}  To describe a semigroup of this type, we need a slightly different  construction, given below.

Suppose that $n=lp(m-1)+p$, $l\geq 2$. Let $X=X_1\cup\dots\cup X_m$ be such a decomposition of $X$ that $|X_1|=p$, $|X_2|=\cdots =|X_m|=lp$. We assume that each of the blocks $X_i$, $2\leq i\leq m$, has, in turn, a decomposition into disjoint subsets  $X_i=U_1^i\cup\dots\cup U_l^i$. Let $G$ be a transitive permutation group acting on $X_2$ of the cardinality $lp$ such that $U_1^i,\dots, U_l^i$ are the imprimitivity blocks, that is, for every block $U_j^2$ we have that $G(U_j^2)$ is again some block $U_t^2$. Let $Z$ be the set of pairs $(a,b)$, where $1\leq a\leq m$, $b=1$ if $a=1$ and $1\leq b\leq l$, if $2\leq a\leq m$. We fix some bijections from $X_1$ to each of $U_i^j$. In this way we identify $X$ with the set $X_1\times Z$, that is the set of triples of the form $(x,a,b)$, where $x\in X_1$, $1\leq a\leq m$, $b=1$ if $a=1$ and $1\leq b\leq l$, if $2\leq a\leq m$.  For $\alpha\in G$ and $(x,2,j)\in X_2$ denote the element $(x,2,j)\alpha$ by $(x_{\alpha},2,j_{\alpha})$. Consider the subsemigroup $T$ of $\IS(Z)$, generated by the following four elements:
$$
\varphi=\bigl( (1,1),(2,1)\bigr]\bigl((2,2)\bigr]\dots \bigl((2,l)\bigr] \bigl((3,1)\bigr]\dots \bigl((3,l)\bigr]\dots \bigl((m,1)\bigr]\dots \bigl((m,l)\bigr],
$$
\begin{multline*}
\psi=\bigl( (1,1)\bigr] \bigl( (2,1), (3,1),\dots, (m,1) \bigr]\bigl( (2,2), (3,2),\dots, (m,2) \bigr]\dots\\ \bigl( (2,l), (3,l),\dots, (m,l) \bigr],
\end{multline*}
$$
g  \text{ -- the idempotent with the domain } \{(1,1)\} \text { and }
$$
$$
h  \text{ -- the idempotent with the domain } Z\setminus\{(1,1)\}.
$$

For $\beta\in T$ and $(a,b)\in \dom(\beta)$ denote $(a,b)\beta$ by $(a_{\beta},b_{\beta})$.
Let $H$ be a subgroup of $G$ stabilizing $U_1^2$. Via the standard restriction and the fixed bijections we can assume that $H$ acts on $X_1$.
The direct product $G\times T$ acts on $X$ by partial permutations as follows:

$$
(x,a,b)\cdot (\alpha,\beta)=\left\lbrace\begin{array}{ll} (x_{\alpha}, a_{\beta}, (b_{\beta})_{\alpha}), & \text{ if } (a,b)\in\dom(\beta) \text{ and } a\neq 1,\\
(x_{\alpha}, 1, 1), & \text{ if } a=b=1, \beta=g \text{ and } \alpha\in H,\\
 (x_{\alpha}, a_{\beta}, (b_{\beta})_{\alpha}), & \text{ if } a=b=1, (a,b)\in\dom(\beta), \beta\neq g\\
 & \text{ and } (b_{\beta})_{\alpha}=b_{\beta}\\
 \text{ not defined},& \text{ otherwise.}
\end{array}\right.
$$

Let $I$ be the ideal in $G\times T$ consisting of the elements with the $T$-coordinate equal to the zero. Similarly as in Construction $(G\times T)/I$ acts on $X$ faithfully, and the cardinality of $(G\times T)/I$ is $2n-p+1$.

{\it Example 2.} Let $n=10$, $p=2$, $l=2$. $X_1=\{1,2\}$, $U_1^2=\{3,4\}$, $U_2^2=\{5,6\}$, $X_2=\{3,4,5,6\}$,  $U_1^3=\{7,8\}$, $U_2^3=\{9,10\}$, $X_3=\{7,8,9,10\}$.
Let $G$ be the four-element cyclic group acting on $X_2$ generated by the cycle $(3,5,4,6)$. The elements of the semigroup $(G\times T)/I$ as the semigroup of partial permutation on $X$ are as follows:
\begin{equation*}
\begin{array}{ll}
\text{\bf the elements of } gSg\setminus\{0\}:&\\
(1)(2)(3](4](5](6](7](8](9](10], & (1,2)(3](4](5](6](7](8](9](10],  \\
\text{\bf the elements of } hSh\setminus\{0\}:&\\
(1](2](3)(4)(5)(6)(7)(8)(9)(10), & (1](2](3,5,4,6)(7,9,8,10), \\
(1](2](3,4)(5,6)(7,8)(9,10), & (1](2](3,6,4,5)(7,10,8,9), \\
(1](2](3,7](4,8](5,9](6,10], & (1](2](3,9](4,10](5,8](6,7], \\
(1](2](3,8](4,7](5,10](6,9], & (1](2](3,10](4,9](5,7](6,8], \\
\text{\bf the elements of } gSh\setminus\{0\}:&\\
(1,3](2,4](5](6](7](8](9](10], & (1,4](2,3](5](6](7](8](9](10], \\
(1,5](2,6](3](4](7](8](9](10], & (1,6](2,5](3](4](7](8](9](10], \\
(1,7](2,8](3](4](5](6](9](10], & (1,8](2,7](3](4](5](6](9](10], \\
(1,9](2,10](3](4](5](6](7](8], & (1,10](2,9](3](4](5](6](7](8] \\
\end{array}
\end{equation*}
and the zero.

{\it \bf Type 5.} The semigroup obtained from the semigroup from Type 4 applying the map $\alpha\mapsto \alpha^{-1}$.

\begin{theorem}
Let $S$ be a semitransitive subsemigroup of $\SIS_n$ of minimal cardinality. Then $|S|=2n-p+1$, where $p$ is the greatest proper divisor of $n$ and, moreover, $S$ is similar to a semigroup constructed in either Type 1, Type 2, Type 3, Type 4 or Type 5.
\end{theorem}
\begin{proof}
Let $S$ be a semitransitive subsemigroup of  $\SIS_n$ of minimal cardinality. From Theorem \ref{th} and the construction above it follows that $|S|=2n-p+1$. Let  $X_1>X_2>\dots >X_m$ be the ordering of the transitivity blocks with respect to the action of $S$ and let $t$ be the cardinality of the smallest of the blocks. From Lemmas~\ref{lemma:lem4} and~\ref{lemma:lem5} and $|(gSg)\setminus\{0\}|+|(hSh)\setminus\{0\}|\geq n$ we have that $t$ is a divisor of $n$, and that $|S|\geq 2n-t+1$. Therefore, it must be $t=p$, $|N|=n-p$ and $|(gSg)\setminus\{0\}|+|(hSh)\setminus\{0\}|=n$.

The restriction $|N|=n-p$ implies that in the notation of the proof of Lemma~\ref{lemma:lem4} the union of all the sets $B_i$, $i\in A$, should have exactly $|B|+|A|-1=m-1$ elements.

The restriction $|(gSg)\setminus\{0\}|+|(hSh)\setminus\{0\}|=n$, together with the result from~\cite{C-O} (recalled in Theorem \ref{th_from_c_o} above), provides us with the information about the structure of $gSg$ and $hSh$. We can therefore assume that for $i\in A$ we have $|X_i|=p$ and for $i\in B$ we also have  $|X_i|=kp$ for some $k\geq 1$.

Let us first show that $k=1$ or $|A|=1$. Let $i\in B$ and $j\in A$. The set $M^{j-i}_{2,1}$, if $j>i$, or the set $M^{i-j}_{1,2}$, if $i>j$, is then non-empty. Then by \eqref{g1} and \eqref{g2} its cardinality should be $kp$. Since $N$ is a disjoint union of exactly $m-1$ different sets $M_{1,2}^j$ or  $M_{1,2}^j$ the cardinality of $N$ must be $(m-1)kp$. From the other hand the cardinality of $N$ is $n-p=tkp+(m-t)p-p$, where $t=|B|$. Hence,
$$
(m-1)kp=tkp+(m-t)p-p.
$$
Therefore, $k=1$, or $m=t+1$. We consider each of these two possibilities separately.

{\it Case 1.} Suppose $k=1$.

{\it Subcase A.} Suppose that there is $i$ such that $i,i+1\in A$.  Let $B_i=\{a_1,\dots, a_{t}\}$. Then it should be $B_{i+1}=\{a_1-1,\dots, a_{t}-1\}$. In addition, from the proof of Lemma \ref{lemma:lem4} and $|N|=n-p$ we should have $|B_i\cup B_{i+1}|=t+1$. It follows  that $B_i=\{a_1,a_1+1,\dots a_1+t-1\}$. Therefore $B=\{j,j+1, j+t-1\}$ for some $j$.

Suppose $|B|>1$. Switching $g$ and $h$ and applying the arguments above we show that $A$ has an analogous structure as $B$. Without loss of generality we can assume that $A=\{1,2,\dots, l\}$, $B=\{l+1,\dots, m\}$.  Since $|A|>1$ and  $|B|>1$ we have $2\leq l\leq m-2$. Assume that $l\geq m-l$ (the other case being similar). Consider $M^{m-l}_{1,2}$. From the restriction $|M^{m-l}_{1,2}|=p$ we obtain that the domain of each element in $M^{m-l}_{1,2}$ contains the blocks $X_{2l-m+1},\dots ,X_l$. Let $x\in X_{2l-m+1}, y\in X_l$, $z\in X_{l+1}$, $t\in X_m$. In $gSh$ there is an element, say $\gamma$, such that it has an arrow from $y$ to $z$, and in $hSh$ there is an element, $\delta$, which has an arrow from $z$ to $t$ and can be decomposed as a product of $m-l-1$ elements. The product $\gamma\delta$ then has an arrow from $y$ to $t$ and belongs to $M^{m-l}_{1,2}\cap gSh\cdot (hSh\cap N_S)$. The domain of such an element must have the empty intersection with $X_{2l-m+1}$. This contradicts the fact that the domain of each element in
$M^{m-l}_{1,2}$ contains the block $X_{2l-m+1}$. This shows  that $|B|=1$.

Suppose $B=\{i\}$, that is, $\dom(h)=X_i$ and $\dom(g)=X\setminus X_i$. Suppose $i\neq 1$ and $i\neq m$. Take $x\in X_{i-1}$, $y\in X_{i}$, $z\in X_{i+1}$ and $\alpha\in gSh$, $\beta\in hSg$ such that $\alpha$ has an arrow from $x$ to $y$, and $\beta$ has an arrow from $y$ to $z$. It follows that $\alpha\beta\in gSg$, and $rank(\alpha\beta)\leq p$. This contradicts to the structure of $hSh$ which we know from Theorem \ref{th_from_c_o}. It follows that  $i= 1$ or $i= m$. Now we know the cardinalities of the transitivity blocks, the positions of $\dom(g)$ and $\dom(h)$ with respect to them, the structures of $gSg$ and $hSh$, also from the proof of Lemma \ref{lemma:lem4} and the limitation $|N|=n-p$ the structure of $N$. Hence, subject to the renumeration of elements in $X$, $S$ must be as in Type 1, or in Type 2.

{\it Subcase B.} Suppose that there is no $i$ such that $i,i+1\in A$. We can assume that there is also no  $i$ such that $i,i+1\in B$, otherwise we switch $g$ and $h$ are are in the previous subcase. Then the transitivity blocks from $\dom(g)$ and $\dom(h)$ are altering, that is, say $A=\{1,3,\dots, 2\lfloor\frac{n-1}{2}\rfloor +1\}$. Again, similarly as in the previous subcase, we conclude that subject to the renumeration of elements in $X$, $S$ must be as in Type 3.

{\it Case 2.} Suppose $m=t+1$ and $k>1$. Let $A=\{i\}$. If $i\neq 1$ and $i\neq m$ we apply similar arguments as in the third paragraph of Subcase A above and obtain a contradiction. Therefore, $i=1$ or $i=m$.

Suppose $i=1$. Consider the elements from $gSh$ with the arrows from $X_1$ to $X_2$. From the restriction on the cardinality of $M^1_{1,2}$  we obtain that there should be $kp$ such elements, each having exactly $p$ arrows from $X_1$ to $X_2$. This implies that there is a decomposition $X_2=U_1\cup\dots\cup U_k$, such that the blocks $U_i$ all have cardinality $p$ and are the intersections of $X_2$ with the images of elements in $M^1_{1,2}$. It follows that for each $i$  it must be that $G(U_i)$ is again some $U_j$, as otherwise $M^1_{1,2}$ would contain more elements than is allowed. This and the other restrictions, similarly as in the two previous subcases, lead to that subject to the renumeration of elements in $X$, $S$ must be as in Type 4. In the case when $i=m$ the arguments are similar, and $S$ must be as in Type 5.
\end{proof}

\noindent{K. Cvetko-Vah, D. Kokol Bukov\v sek, T. Ko\v sir: Department of Mathematics,
University of Ljubljana, Jadranska 19, SI-1000 Ljubljana, Slovenia.\\ e-mail: karin.cvetko@fmf.uni-lj.si,
damjana.kokol@fmf.uni-lj.si,\\ tomaz.kosir@fmf.uni-lj.si.}

\medskip

\noindent{G. Kudryavtseva: Centre for systems and information technologies,
University of Nova Gorica, Vipavska cesta 13, p.p. 301, SI-5001 Nova Gorica, Slovenia. \\ e-mail:  ganna.kudryavtseva@p-ng.si.}


\begin{thebibliography}{99}

\bibitem{BDKKO} J.~Bernik, R. Drnov\v sek, D. Kokol Bukov\v sek,
T. Ko\v sir, and M. Omladi\v c.  \emph{Reducibility and triangularizability of semitransitive spaces of operators}. Houston J. Math. {\bf 34} (2008), no. 1, 235--247.

\bibitem{BDKKOR} J.~Bernik, R. Drnov\v sek, D. Kokol Bukov\v sek,
T. Ko\v sir, M. Omladi\v c, and H.~Radjavi.
\emph{On semitransitive Jordan algebras of matrices}. Preprint.

\bibitem{BGMRT} J.~Bernik, L.~Grunenfelder, M.~Mastnak, H.~Radjavi, and V.~G.~Troitsky.
\emph{On semitransitive collections of operators}. Semigroup Forum {\bf 70} (2005), 436--450.



\bibitem{C-O} K. Cvetko-Vah, D. Kokol Bukov\v{s}ek, T. Ko\v{s}ir, G. Kudryavtseva, Ya.~Lavrenyuk and A.
Oliynyk, Semitransitive subsemigroups of the symmetric inverse
semigroups. {\it Semigroup Forum}, {\bf 78} (2009),138--147.

\bibitem{East} J. East, A presentation of the singular part of the symmetric inverse monoid. {\it Comm. Algebra} {\bf 34} (2006), no. 5, 1671--1689.

\bibitem{GM} O. Ganyushkin, V. Mazorchuk, Classical finite transformation semigroups, an introduction. Algebra and Applications, Vol. {\bf 9}, Springer Verlag, 2009.

\bibitem{MM} V. Maltcev, V. Mazorchuk,  Presentation of the singular part of the Brauer monoid. {\it Math. Bohem.} {\bf 132} (2007), no. 3, 297--323.

\bibitem{RaT} H.~Radjavi and V.~G.~Troitsky. \emph{Semitransitive subspaces of operators}. To appear
in Linear and Multilinear Algebra.

\bibitem{RT} H. Rosenthal and V. G. Troitsky. \emph{Strictly semi-transitive operator algebras}.
Journal of Operator Theory  {\bf 53}  (2005), 315--329.

\bibitem{Bled}
Semitransitivity Working Group at LAW'05, Bled. \emph{Semitransitive subspaces of matrices}.
Electronic Journal of Linear Algebra {\bf 15} (2006), 225-238.
\end{thebibliography}
\end{document}